\pdfoutput=1
\documentclass[11pt, oneside]{amsart}

\usepackage{hyperref}
\hypersetup{%
	pdfpagemode={UseOutlines},
	bookmarksopen,
	pdfstartview={FitH},
	colorlinks,
	linkcolor={blue},
	citecolor={blue},
	urlcolor={blue}
}

\usepackage{geometry} 

\usepackage{dcolumn}
\newcolumntype{d}{D{.}{.}{-1} } 

\usepackage{amsmath, amssymb, amsthm}
\usepackage{mathrsfs}
\usepackage{enumitem}
\usepackage{aligned-overset}
\usepackage{mathtools}
\usepackage[numbers]{natbib}

\numberwithin{equation}{section}

\newcommand\eps\varepsilon
\newcommand\newk{\textsc{k}}
\newcommand\newj{\textsc{j}}
\newcommand\newalpha{{\widehat{\alpha}}}
\newcommand\newg{\widehat{g}}
\newcommand\UU{{\mathcal U}}
\renewcommand\Re{\operatorname{Re}}
\renewcommand\Im{\operatorname{Im}}
 \DeclareMathOperator{\Tr}{tr}
 
\newtheorem{theorem}{Theorem}[section]

\newtheorem{lemma}[theorem]{Lemma}
\newtheorem{corollary}[theorem]{Corollary}
\newtheorem{proposition}[theorem]{Proposition}

\theoremstyle{remark}
\newtheorem{remark}[theorem]{Remark}

\newcommand{\R}{\mathbb{R}} 
\newcommand{\N}{\mathbb{N}} 
\newcommand{\V}{\mathcal{V}} 
\newcommand{\Fock}{\mathcal{F}} 
\newcommand{\C}{\mathbb{C}} 

\title[The isoperimetric inequality for partial sums of Toeplitz eigenvalues]{The isoperimetric inequality for partial sums of Toeplitz eigenvalues in the Fock space}
\author{Fabio Nicola}
\address[Fabio Nicola]{Dipartimento di Scienze Matematiche, Politecnico di Torino, Corso Duca degli Abruzzi 24, 10129 Torino, Italy}
\email{fabio.nicola@polito.it}

\author{Federico Riccardi}
\address[Federico Riccardi]{Dipartimento di Scienze Matematiche, Politecnico di Torino, Corso Duca degli Abruzzi 24, 10129 Torino, Italy}
\email{federico.riccardi@polito.it}

\author{Paolo Tilli}
\address[Paolo Tilli]{Dipartimento di Scienze Matematiche, Politecnico di Torino, Corso Duca degli Abruzzi 24, 10129 Torino, Italy}
\email{paolo.tilli@polito.it}
\begin{document}

    
    \begin{abstract}
        \noindent We prove that, among all subsets $\Omega\subset \C$ having
        circular symmetry and prescribed measure, the ball is the only maximizer of the sum of the first $\newk$ eigenvalues ($\newk\geq 1$) of the corresponding Toeplitz operator $T_\Omega$ on the Fock space $\Fock$. As a byproduct, we prove that balls maximize any Schatten $p$-norm of $T_\Omega$ for $p>1$ (and minimize 
        the corresponding quasinorm for $p<1$), and that the second eigenvalue is maximized by a particular annulus. Moreover, we extend some of these results to general radial symbols in $L^p(\C)$, with $p > 1$,  characterizing those that maximize the sum of the first $\newk$ eigenvalues. We also show a symmetry breaking phenomenon for the second eigenvalue, when the assumption of circular symmetry is dropped. 
    \end{abstract}
    
    \maketitle
    
    \section{Introduction}

    Recent years have seen  a growing interest in optimization problems arising in the context of time-frequency analysis and complex analysis, for example concerning the maximal concentration of the short-time Fourier transform with Gaussian window \cite{gomez_guerra_ramos_tilli, trapasso, nicolatilli_fk} or the wavelet transform with respect to the Cauchy wavelet \cite{ramos-tilli} and related localization operators \cite{galbis2022norm, nicola_riccardi_quantitative, nicolatilli_norm, riccardi_optimal_estimate}, optimal estimates for functionals and contractive estimates on various reproducing kernel Hilbert spaces \cite{garcia_ortega_stability, frank2023sharp, frank2023generalized,  kalaj1, kulikov,ortega,nicola_stabilizer,NRT2024}. A case of special interest is given by the Fock space $\Fock$ of entire functions $F:\C\to\C$ with finite norm 
    \[
\|F\|_{\Fock}^2=\int_\C |F(z)|^2 e^{-\pi |z|^2} dA(z),
    \]
   where $dA(z)$ stands for the Lebesgue measure in the plane. In \cite{nicolatilli_fk} it was proved that, for every 
    normalized $F\in \Fock$  and every set $\Omega \subset \C$ of finite two-dimensional Lebesgue measure  (denoted by $|\Omega|$, and simply called area or measure of $\Omega$), the concentration inequality
	\begin{equation}
    \label{concineq1}
		\int_{\Omega} |F(z)|^2 e^{-\pi |z|^2}\, dA(z) \leq 1 - e^{-|\Omega|} 
	\end{equation}
 holds true, and equality is achieved if and only if $\Omega$ is (up to a negligible set) a ball
 centered at some $w\in\C$
 and $F$ is a normalized multiple of
 the reproducing kernel $e^{\pi\overline{w} z}$. This result can be seen as the optimal norm estimate
 $
 \Vert T_\Omega \Vert_{\Fock \rightarrow \Fock} \leq 1 - e^{-|\Omega|}
 $
 for the Toeplitz operator $ T_{\Omega} = P \chi_{\Omega} \colon \Fock \rightarrow \Fock$, where $\chi_{\Omega}$ is the characteristic function of the set $\Omega$ and $P$ is the projection operator from $L^2(\C, e^{-\pi |z|^2} dA(z))$ onto $\Fock$.
Indeed 
$T_\Omega$, also known as \emph{time-frequency localization operator} \cite{daubechies}, is self-adjoint, compact and positive on $\Fock$ (see e.g. \cite{zhu_book}), and hence, labeling its eigenvalues in a decreasing way
as 
\begin{equation}
\label{specT}
\lambda_1(\Omega) \geq \lambda_2(\Omega) \geq  \cdots > 0,
\end{equation}
and interpreting
the left-hand side of 
\eqref{concineq1} as a Rayleigh quotient, one obtains
    \begin{equation}
    \label{fk2}
        \lambda_1(\Omega) \leq 1 - e^{-|\Omega|}\qquad\text{(with equality $\iff$ $\Omega$ is a ball),}
    \end{equation}
that is a ``Faber--Krahn inequality''
for the operator $T_\Omega$: 
among all sets $\Omega$ of prescribed measure, 
the first eigenvalue $\lambda_1(\Omega)$ of $T_\Omega$ is maximized when $\Omega$ is a ball.

\begin{remark}\label{rem1}
    The classical Faber--Krahn inequality \cite{henrot} states that, among sets of given measure, the ball \emph{minimizes} 
    the first eigenvalue $\lambda_1^D(\Omega)$ of the Dirichlet Laplacian, while in \eqref{fk2} the ball \emph{maximizes}
    $\lambda_1(\Omega)$. As explained in Remark 1.3 of \cite{nicolatilli_fk}, however, 
    the two Faber--Krahn inequalities go in the same direction: indeed, the spectrum \eqref{specT}
    of $T_\Omega$ (where $\lambda_n(\Omega)\to 0$ and  $\lambda_1(\Omega)$ is the \emph{largest} eigenvalue)
    should not be compared
    with the spectrum of the Dirichlet Laplacian $-\triangle$ (where 
    $\lambda_n^D(\Omega)\to \infty$ and $\lambda_1^D(\Omega)$ is the \emph{smallest} eigenvalue) but, rather, with the spectrum of its \emph{inverse} $(-\triangle)^{-1}$ as a compact operator on $H^1_0(\Omega)$,
  with \emph{decreasing}  eigenvalues $1/\lambda_n^D(\Omega)\to 0$. This parallelism between the eigenvalues $\lambda_n(\Omega)$ of $T_\Omega$ and the inverse eigenvalues $1/\lambda_n^D(\Omega)\to 0$ of the Laplacian
will be pursued further in the following.
\end{remark}

The study of the spectrum of localization operators was initiated in
\cite{daubechies}, and several results are now available concerning the eigenvalue distributions and other global spectral properties
(see, e.g., \cite{abreu,kulikov2024} and the references therein). However, contrary to the case of the Laplacian 
(where,  since the pioneering work of P\'olya and Szeg\H{o} \cite{polya_szego}, 
 spectral optimization problems
have  become a well-established research field, see \cite{ashbaugh_benguria,grebenkov_nguyen,henrot}
for some modern accounts), much less is known about extremal properties 
of the eigenvalues $\lambda_n(\Omega)$, 
and the Faber--Krahn inequality in \cite{nicolatilli_fk} was the
first step in this direction 
(prior to \cite{nicolatilli_fk}, this had been proved in \cite{galbis2022norm} 
but only for sets with circular symmetry).

The goal
of this paper is to carry over the program initiated in \cite{nicolatilli_fk} with the Faber--Krahn
inequality \eqref{fk2}, by investigating (in the spirit of \cite{polya_szego})
further basic questions, such as:  among all sets $\Omega\subset\C$ of \emph{prescribed area}, which sets  maximize the second
eigenvalue $\lambda_2(\Omega)$? Which sets maximize the sum 
\begin{equation}
    \label{sommaautov}
\sum_{n=1}^\newk \lambda_n(\Omega)
\end{equation}
of the first $\newk$ eigenvalues? Or, again, which sets
maximize the Schatten norms
\begin{equation}
\label{defschatten}
  \left(\sum_{k=1}^\infty \lambda_k(\Omega)^p\right)^{\frac 1 p}\quad 
\text{?}
\end{equation}
One reason of interest is that 
the machinery developed in \cite{nicolatilli_fk} for Faber--Krahn 
(though adaptable to contexts other than the Fock space, see,
e.g., \cite{frank2023sharp,frank2023generalized, kulikov, ortega,ramos-tilli}) 
is \emph{not suitable} to attack the aforementioned problems, 
since the differential inequality obtained in 
\cite{nicolatilli_fk} is no longer available: therefore,
new ideas and tools must be introduced. Moreover, contrary to \eqref{fk2}, for the maximization of $\lambda_2(\Omega)$
circular symmetry \emph{breaks down}, in that the optimal set $\Omega$   (if one exists) is no longer rotationally invariant around any point:
this symmetry breaking suggests
that it might be extremely difficult to characterize 
the sets $\Omega$ of prescribed area  that maximize $\lambda_2(\Omega)$
(or the sum of the first $\newk$ eigenvalues, or the Schatten norms etc.)
without additional constraints, and 
these questions remain (in their full generality) three 
challenging open problems. 
We are able, however, to solve
these problems within the class of  sets (of given area) 
that are \emph{circularly symmetric}, i.e. measurable sets $\Omega\subset\C$  that, possibly after a translation, have the form
\begin{equation}
    \label{circsym}
    \Omega=\left\{z\in\C\quad\text{such that}\,\,\, |z|\in E\right\},\quad\text{for some  measurable set $E\subset\R^+$.}
\end{equation}
Any such set $\Omega$ is obtained by rotating the set $E\subset \R^+$
(thought of as a subset of the positive real halfline) around the origin of the complex plane $\C$, and the area constraint $|\Omega|=s$ corresponds
to the integral condition
\[
2\pi \int_E r \,dr=s.
\]
When $\Omega$ is \emph{radial}, that is, circularly symmetric \emph{around the origin} as in \eqref{circsym}, it is well known \cite{daubechies, seip}
that the eigenfunctions of $T_\Omega$ are the normalized 
monomials 
	\begin{equation}
    \label{monomials}
		e_k(z) = \dfrac{\pi^{k/2}\, z^k}{\sqrt{k!\,}},  
        \quad k=0,1,2,\ldots
	\end{equation}
(which also form an orthonormal basis of the Fock space $\Fock$),
and the corresponding eigenvalues are given by the Rayleigh quotients
\begin{equation}
\label{eq:expression eigenvalues}
\langle T_\Omega \,e_k,e_k\rangle_\Fock=        \int_{\Omega} |e_k(z)|^2 e^{-\pi|z|^2} \, dA(z), \quad k=0,1,2,\ldots
\end{equation}
Of course, there is a one-to-one correspondence
between these numbers and the eigenvalues labeled decreasingly as in
\eqref{specT},
but which $k$ corresponds to which $n$
strongly depends on the set $\Omega$: for instance, if $\Omega$ is a thin
annulus of radius $R\gg 1$, then $\lambda_1(\Omega)$ is likely to
correspond to some $k\sim \pi R^2$, since $|e_j(z)|^2$ has its peaks
along a circle of radius $\sqrt{j/\pi}$. If $\Omega$ is a ball,
however, one has (cf. \cite{daubechies,seip})
\begin{equation}\label{eigenval_ball0}
\lambda_{n}(\Omega)= 
    \int_{\Omega} |e_{n-1}(z)|^2 e^{-\pi|z|^2} \, dA(z)=\int_0^{\pi r^2} \frac{s^{n-1}}{(n-1)!} e^{-s}\, ds, \qquad n\geq1,  
    \end{equation}
where $r$ is the radius of the ball, and in particular (cf. \eqref{fk2})
\begin{equation}
    \label{eigenval_ball}
\lambda_1(B)=1-e^{-|B|},\quad\text{for every ball $B$ of area $|B|$.}
\end{equation}
The origin as center of symmetry is not privileged: the case where $\Omega$ has circular symmetry around an \emph{arbitrary}
point $w\in\C$ can be easily reduced to the case where $w=0$, thanks
to the following general fact.
\begin{remark}[Translation invariance]\label{reminv} If 
a measurable set $\Omega\subset\C$ is translated
by a vector $w\in \C$, then
the spectrum of $T_\Omega$ is unchanged, i.e. $\lambda_n(\Omega)
=\lambda_n(\Omega+w)$ for every $n\geq 1$. This follows from
the well-known fact (see e.g. \cite{zhu_book}) that the operator $\UU_{w}:\Fock\to\Fock$ 
defined by
\begin{equation}
\label{unitaryop}
			\UU_{w}F(z) := e^{-\,\frac{ \pi|w|^2}2+ \pi z \overline{w}}
            F(z-w), \quad F \in \Fock,\ z \in \C,
\end{equation}
is \emph{unitary} (in fact, $|\UU_w F(z)|^2 e^{-\pi |z|^2}=
|F(z-w)|^2 e^{-\pi|z-w|^2}$), so that the quadratic forms induced by $T_{\Omega}$ and $T_{\Omega+w}$ are related by
$
\langle T_\Omega F,\,F\rangle_\Fock = \langle T_{\Omega+w} \,\UU_wF,\,\UU_wF\rangle_\Fock
$. Thus, if $F_n$ is an eigenfunction of 
$T_\Omega$ relative to $\lambda_n(\Omega)$, then $\UU_wF_n$ is
the corresponding eigenfunction for $T_{\Omega+w}$; in particular,
if a set $\widehat{\Omega}$ has circular symmetry around an
\emph{arbitrary} point
$w\in\C$ (i.e. if $\widehat\Omega=\Omega+w$ for some $\Omega$
as in \eqref{circsym}), then the eigenvalues $\lambda_n(\widehat\Omega)$
are given by \eqref{eq:expression eigenvalues}, while the eigenfunctions
of $T_{\widehat\Omega}$ are the functions $\UU_w e_k(z)$, cf. \eqref{monomials}.
\end{remark}
After these preliminaries, our main results  can be
stated as follows. Concerning the maximization of $\lambda_2(\Omega)$,
among sets with rotational symmetry the optimal set 
is always
a certain \emph{annulus}, whose shape and size depend on 
the prescribed area: 
\begin{theorem}[Krahn--Szeg\H{o} inequality for circularly symmetric sets]\label{teolambda2} Let $\Omega\subset\C$ be a circularly symmetric
set of finite area $s:=|\Omega|>0$. Then 
\begin{equation}
    \label{inequality lambda 2}
\lambda_2(\Omega)\leq \lambda_2(A_s)=  e^{s/(1-e^s)}\left(1-e^{-s}\right),
\end{equation}
where $A_s$ is
the annulus
$\left\{ r < |z| < R\right\}$ of area $s$,
with radii  defined by the equations
\begin{equation}
    \label{defradii}
\pi r^2=\displaystyle\frac{s}{e^s-1},\quad
\pi R^2=\displaystyle\frac{s e^s}{e^s-1}.
\end{equation}
Moreover, equality in \eqref{inequality lambda 2} occurs if and only if $\Omega$ is equivalent, up to a translation, to
the annulus $A_s$.
\end{theorem}
Here and in the following we say that two measurable subsets $\Omega,\Omega'\subset\C$ are equivalent if their symmetric difference $\Omega\Delta\Omega':=(\Omega\setminus \Omega') \cup (\Omega'\setminus\Omega)$ has measure zero. 

The name ``Krahn--Szeg\H{o}'' is by formal
analogy, along the lines of Remark \ref{rem1},  with
the minimization problem for
the second eigenvalue of the
Laplacian (\cite{ashbaugh_benguria,henrot}).
As we will show, the optimal annulus $A_s$ defined by \eqref{defradii} arises, in quite a natural way, as a superlevel set
\begin{equation}
\label{defAs}
A_s=\left\{ z\in\C\text{\ such that\ }\,\, |e_1(z)|^2 e^{-\pi |z|^2} > t\right\}
\end{equation}
where $e_1(z)$ is as in \eqref{eq:expression eigenvalues}, and $t$ chosen as to match the prescribed area $|A_s|=s$. For these annuli
one has $\lambda_2(A_s)=\lambda_1(A_s)$, and the corresponding 
eigenfunction is any linear combination of $e_0(z)$ and $e_1(z)$.
As we mentioned, however, these annuli are
no longer optimal (at least for small area $s$) among planar sets
without circular symmetry:
\begin{proposition}[Symmetry breaking for $\lambda_2(\Omega)$]\label{prononrad}
    For every small enough $s>0$, there exist sets $\Omega\subset\C$ of area $s$ (such as the union
    of two equal disks sufficiently apart from each other) that violate \eqref{inequality lambda 2}, i.e. such that
    \begin{equation}
\label{tesiprop}        
    \lambda_2(\Omega)>  e^{s/(1-e^s)}\left(1-e^{-s}\right).
    \end{equation}
\end{proposition}
It is an open problem to establish whether \eqref{inequality lambda 2} can
also be violated by sets $\Omega$ of large measure $s$ (if not, then $A_s$ would
maximize $\lambda_2(A_s)$ without any symmetry constraint). Moreover,
without symmetry constraints, it
is not even known (no matter how $s>0$ is chosen) whether a set $\Omega$
exists, that maximizes $\lambda_2(\Omega)$ among sets of area $s$
(the union of two disjoint disks, which by the way optimizes the Krahn--Szeg\H{o} inequality for the Laplacian, satisfies \eqref{tesiprop}
for small $s$, but does \emph{not} maximize $\lambda_2(\Omega)$ among
sets of area $s$, see the proof of Proposition \ref{prononrad}).

Summing up, among circular sets of a area $s$, $\lambda_1(\Omega)$ and $\lambda_2(\Omega)$
are maximized, respectively, by the ball (which is also optimal among arbitrary sets of area $s$) and by a the annulus $A_s$.
The sum $\lambda_1(\Omega)+\lambda_2(\Omega)$, however,
is maximized again by the
ball and, more generally, the same is true for the sum in \eqref{sommaautov}:
\begin{theorem}[Sum of the first $\newk$ eigenvalues for circularly symmetric sets]\label{teoPS} 
Let $\Omega\subset\C$ be a circularly symmetric
set of finite area $|\Omega|>0$. Then, for every integer $\newk\geq 1$,
\begin{equation}
\label{ineqPS}
\sum_{k=1}^{\newk} 
\lambda_k(\Omega) \leq
\sum_{k=1}^{\newk} 
\lambda_k(\Omega^*)
=\newk - \sum_{k=0}^{\newk-1} (\newk-k)\dfrac{|\Omega|^k}{k!}e^{-|\Omega|},
\end{equation}
where $\Omega^*$ is the ball such that $|\Omega^*|=|\Omega|$. 
Moreover, equality holds if and only if $\Omega$ is equivalent
(up to a translation)
to $\Omega^*$.
\end{theorem}
As is well known, the sum $\sum_{k=1}^{\newk} 
\lambda_k(\Omega)$ is the trace of the operator $T_\Omega$
\emph{restricted} to the space generated by its first $\newk$ eigenfunctions;
more precisely, if $P_V$ is the orthogonal
projector from $\Fock$ to any $\newk$-dimensional subspace $V\subset\Fock$
spanned by $\newk$ orthogonal functions $F_1,\ldots,F_\newk$, one has
\begin{equation}
\label{minmax1}
 \Tr (P_V T_\Omega P_V) = \sum_{k=1}^\newk \langle
T_\Omega F_k,F_k \rangle =
\sum_{k=1}^\newk  \int_{\Omega} |F_k(z)|^2 e^{-\pi |z|^2} \, dA(z) \leq
\sum_{k=1}^{\newk} 
\lambda_k(\Omega),
\end{equation}
with equality when $V$ is spanned by the \emph{first} $\newk$ eigenfunctions of
$T_\Omega$. 
In fact, we obtain \eqref{ineqPS} as a consequence of the following more precise statement, which also takes into account the center of symmetry of $\Omega$
and (implicitly) the first $\newk$ eigenfunctions of $T_\Omega$:
\begin{theorem}\label{th:main}
		Let $\Omega \subset \C$ be a set with finite measure
        $|\Omega|>0$, having circular symmetry around some point $w\in\C$.
If $F_1, \ldots, F_\newk$ are any $\newk$
         orthonormal functions in the Fock space $\Fock$,
        then
		\begin{equation}\label{eq:main estimate}
			\sum_{k=1}^\newk  \int_{\Omega} |F_k(z)|^2 e^{-\pi |z|^2} \, dA(z) \leq 
            \newk - \sum_{k=0}^{\newk-1} (\newk-k)\dfrac{|\Omega|^k}{k!}e^{-|\Omega|}.
        \end{equation}
        Moreover, equality is achieved if and only if $\Omega$ is equivalent to a ball centered at $w\in\C$ and,
        for some unitary matrix $U \in \C^{\newk \times \newk}$, there holds
		\begin{equation}\label{unitU}
			\begin{pmatrix}
				F_1(z) \\
				\vdots \\
				F_\newk(z)
			\end{pmatrix} = U \begin{pmatrix}
				\UU_w e_0(z) \\
				\vdots \\
				\UU_w e_{\newk-1}(z)
			\end{pmatrix}\quad \forall z\in\C,
		\end{equation}
where $\UU_w:\Fock\to\Fock$ is the unitary operator  defined in
\eqref{unitaryop} and the $e_k$'s are the monomials defined
in \eqref{monomials}.
	\end{theorem}
The validity of \eqref{ineqPS} and \eqref{eq:main estimate} without
the assumption that $\Omega$ has circular symmetry remains
an \emph{open problem}; in fact,
we conjecture that here 
 (contrary to Proposition \ref{prononrad})
 no symmetry breaking occurs, i.e. that the ball maximizes
 the sum $\sum_{k=1}^\newk \lambda_k(\Omega)$ among
 \emph{all sets} of given area (this would extend the Faber--Krahn
 inequality of \cite{nicolatilli_fk}, from $\lambda_1(\Omega)$
 to the sum $\sum_{k=1}^\newk \lambda_k(\Omega)$).
 With this respect, let us mention that also the corresponding
 problem for the eigenvalues of the Laplacian, i.e. the maximization
 of $\sum_{k=1}^\newk \lambda_k^D(\Omega)^{-1}$  under
 an area constraint (cf. Remark \ref{rem1}) is still open
 when $\newk>1$
 (see \cite{henrot}; if, instead of the area constraint, one
 considers a conformal radius constraint on $\Omega$, then the ball \emph{minimizes} the same sum: this is the P\'olya--Schiffer inequality, see \cite{henrot,polya_schiffer}). When $\newk=\infty$, 
 Luttinger \cite{luttinger} proved that, under an area constraint,
 the ball maximizes
 $\sum_{k=1}^\infty \lambda_k^D(\Omega)^{-1}$, 
 but the same
 problem for localization operators is trivial, since the sum
 $\sum_{k=1}^\infty \lambda_k(\Omega)$, i.e. the trace of $T_\Omega$,
 is equal to the area $|\Omega|$ for \emph{any} measurable set $\Omega$ of finite area (see \cite{zhu_book}).
 The problem becomes nontrivial (even with the circular symmetry
 constraint) if one considers infinite sums as in \eqref{defschatten}
 with an exponent $p\not=1$ or, more generally, the trace
 $\Tr \Phi(T_\Omega)$ where where $\Phi$ is an arbitrary convex function:
\begin{theorem}[Trace inequality for circularly symmetric sets]\label{teotrace}
   Let $\Omega\subset\C$ be a circularly symmetric
set of finite area $|\Omega|>0$, and let   $\Phi \colon (0,1) \to \R$ 
be a convex function. Then $\Tr \Phi(T_\Omega)\leq \Tr \Phi(T_{\Omega^*})$,
that is,
\begin{equation}
\label{eq:inequality traces}
    \sum_{k=1}^{\infty} \Phi(\lambda_k(\Omega))
    \leq\sum_{k=1}^{\infty} \Phi(\lambda_k(\Omega^*)),
\end{equation}
where $\Omega^*$ is the ball such that $|\Omega^*|=|\Omega|$.
Moreover, if $\Tr \Phi(T_{\Omega^*})$ is finite and if $\Phi$ is not affine in the interval $(0,\lambda_1(\Omega^*))$, then equality in \eqref{eq:inequality traces} is achieved if and only if $\Omega$ is equivalent, up to a translation, to $\Omega^*$.
 \end{theorem}
 In particular, letting $\Phi(t)=t^p$ when $p>1$, or $\Phi(t)=-t^p$ when $p\in (0,1)$, one obtains the following result: 
 \begin{corollary}
     Among circularly symmetric sets $\Omega$ of given area, the ball
     maximizes all the Schatten norms in \eqref{defschatten} when $p>1$,
     and minimizes all the Schatten quasi-norms in \eqref{defschatten} when $p\in (0,1)$.
 \end{corollary}
Again, the validity of these results without the assumption of circular symmetry
is an open problem, except when $p=2$ (Hilbert--Schmidt norm): indeed,
       in \cite{nicola_riccardi_quantitative}
       it was proved that balls are maximizers,
       without the assumption that $\Omega$ is circularly symmetric. 
       
We point out
       for future reference that, as already observed, for every set $\Omega$ of finite measure,
       there holds
\begin{equation}
    \label{traccia}
    \sum_{k=1}^\infty \lambda_k(\Omega)=|\Omega|<\infty
\end{equation}       
(in particular $T_\Omega$ is of trace class, so that all the above Schatten
norms are finite, see \cite{zhu_book}).
 
Finally we observe that, probably, a quantitative form of the above results could be established by refining these methods. Further, similar inquiries can be explored in higher dimensions or alternative frameworks, such as functions within weighted Bergman spaces or polynomial spaces. We intend to address these questions in future research.

\section{Proof of Theorems \ref{teoPS} and 
\ref{th:main}}\label{sec:proof of main theorem}

As explained in the Introduction, Theorem \ref{teoPS} follows
from Theorem \ref{th:main}, whose proof requires
some preliminary considerations and a lemma which might be of 
some interest in other contexts as well.

First, by Remark \ref{reminv}, it suffices to prove Theorem \ref{th:main}
when $w=0$, i.e. when $\Omega$ is circularly symmetric around the
origin: the general case then follows by replacing the 
$\newk$ orthogonal functions $\{F_j\}$ with $\{\UU_w F_j\}$ which,
as explained after \eqref{unitaryop}, corresponds to
a shift $\Omega\mapsto \Omega-w$, 
i.e. to the change of variable $z \to z-w$ in all
the integrals in \eqref{eq:main estimate}.

Now, given $\Omega$ and $F_1,\ldots,F_{\newk}$ as in Theorem \ref{teoPS}
with $w=0$, our goal is to prove that
\begin{equation}
    \label{goal1}
\sum_{k=1}^\newk  \int_{\Omega} |F_k(z)|^2 e^{-\pi |z|^2} \, dA(z) \leq
\sum_{k=0}^{\newk-1}  \int_{\Omega^*} |e_k(z)|^2 e^{-\pi |z|^2} \, dA(z),
\end{equation}
where $\Omega^*$ is the ball (centered at the origin) such that
$|\Omega^*|=|\Omega|$, and the $e_k$'s are the first $\newk$ monomials 
in \eqref{monomials}. Since $\Omega^*$ is a ball, due to \eqref{eigenval_ball0} the second sum in \eqref{goal1} coincides
with $\sum_{k=1}^\newk \lambda_k(\Omega^*)$ which in turn, by a direct 
computations of the integrals in \eqref{eigenval_ball0}, coincides
with the right-hand side of \eqref{eq:main estimate}: in other
words, \eqref{goal1}
is equivalent to \eqref{eq:main estimate}. To complete the proof
of Theorem \ref{th:main}, we will also prove that equality in \eqref{goal1}
occurs only when $\Omega=\Omega^*$ and
\begin{equation}
\label{unitU2}
 			\begin{pmatrix}
 				F_1(z) \\
 				\vdots \\
 				F_\newk(z)
 			\end{pmatrix} = U \begin{pmatrix}
 				e_0(z) \\
 				\vdots \\
 				e_{\newk-1}(z)
 			\end{pmatrix}
 		\end{equation}
with $U$ unitary matrix, which corresponds to \eqref{unitU} when $w=0$.

Since $\Omega$ is circularly symmetric around the origin,
every eigenvalue $\lambda_n(\Omega)$ can be found in \eqref{eq:expression eigenvalues}  corresponding to some value of $k$
(the correspondence $n\mapsto k_n$, as already mentioned,
depends on $\Omega$), and so, if we use
\eqref{minmax1} to bound the left-hand side of \eqref{goal1},
we obtain that
\[
\sum_{k=1}^\newk  \int_{\Omega} |F_k(z)|^2 e^{-\pi |z|^2} \, dA(z) \leq
\sum_{n=1}^{\newk}  \int_{\Omega} |e_{k_n}(z)|^2 e^{-\pi |z|^2} \, dA(z),
\]
for some unknown integers $k_n\geq 0$, which we may sort in such a way that
\begin{equation}
    \label{sort1}
0\leq k_1 < k_2 <\cdots < k_\newk.
\end{equation}
Thus, to prove \eqref{goal1}, it suffices to prove that
\begin{equation}
    \label{goal2}
    \sum_{n=1}^{\newk}  \int_{\Omega} |e_{k_n}(z)|^2 e^{-\pi |z|^2} \, dA(z)
    \leq
    \sum_{k=0}^{\newk-1}  \int_{\Omega^*} |e_k(z)|^2 e^{-\pi |z|^2} \, dA(z),
\end{equation}
for every $\newk$-tuple of integers $k_n$ satisfying \eqref{sort1}.
To simplify the notation in view of \eqref{goal2}, given a vector of $\newk$ integers
$\alpha=(k_1,\ldots,k_\newk)$ with the $k_i$'s as in \eqref{sort1}
we let
\begin{equation}
\label{not1}
    u_\alpha(z):=e^{-\pi |z|^2}\sum_{n=1}^{\newk}  |e_{k_n}(z)|^2,
    \quad u_*(z):= e^{-\pi |z|^2}\sum_{k=0}^{\newk-1}  |e_{k}(z)|^2,
\end{equation}
where $u_*(z)$ is just the particular case of $u_\alpha$
when $\alpha=(0,1,\ldots,\newk-1)$, as in the right-hand side of \eqref{goal2}.
\begin{remark}\label{remsuplevel} Since $u_\alpha(z)$ is
a continuous radial function, any of its \emph{superlevel sets}
\begin{equation}
    \label{suplev}
E_t =\left\{z\in \C\,|\,\, u_\alpha(z) > t\right\},\quad t> 0,
\end{equation}
is an open set with circular symmetry, i.e. the union of 
a ball and countably many annuli (the ball and some annuli may
be empty); in fact, since along every ray $u_\alpha(z)$  has finitely
many critical points, the annuli are in finite number, 
so that every superlevel set has the form
\begin{equation}
\label{structure}
    \left\{|z| < R_0\right\}\cup \bigcup_{i=1}^m 
    \left\{r_i< |z|<R_i\right\},\quad 0\leq R_0 \leq r_1 \leq R_1
    \leq \cdots \leq r_m\leq R_m
\end{equation}
for some $m\geq 1$ (the cases where $R_0=0$ or $r_i=R_i$ allow
for empty sets). In view of \eqref{goal2}, superlevel sets are relevant
in that, for any set $\Omega$ of finite measure $|\Omega|>0$, there holds
\begin{equation}\label{bathtub}
    \int_\Omega u_\alpha(z)\,dA(z)\leq
    \int_E u_\alpha(z)\,dA(z)\quad\text{($E$ superlevelset of
    $u_\alpha$ such that $|E|=|\Omega|$),}
\end{equation}
i.e. any integral increases on passing from a set $\Omega$ to a superlevel
set with the same measure (this well known property --valid in general
and not peculiar of $u_\alpha$--
is a particular case of the so called ``bathtub principle'', see
\cite{liebloss}). Let us just observe that, since 
$u_\alpha$ has no flat zones, its
distribution
function
$t\mapsto \left|E_t\right|$
(which maps the interval $(0,\max u_\alpha)$ surjectively onto 
$(0,+\infty)$) is continuous: therefore, for  any set
$\Omega$ of finite measure, one can always find a $t>0$ 
such that $\left\vert E_t\right\vert=|\Omega|$, i.e. a superlevel
set $E$ of $u_\alpha$ always exists, that can be used
in \eqref{bathtub} matching the measure of $\Omega$
(it is easy to see that such $E$ is unique, and that the inequality
in \eqref{bathtub} is strict unless $|E\triangle\Omega|=0$).
\end{remark}
Using the bathtub principle  \eqref{bathtub}
(and the notation in \eqref{not1} for $u_*$) we see that
\eqref{goal2} follows from
\begin{equation}
\label{goal3}
\int_E u_\alpha(z)\,dA(z)\leq \int_{\Omega^*} u_*(z)\,dA(z)
\quad\text{($E$ superlevel set of $u_\alpha$ such that
$|E|=|\Omega^*|$),}
\end{equation}
which we will prove (on passing from \eqref{goal2} to \eqref{goal3}
we loose track of $\Omega$ except for its measure,
implicitly retained
through the ball $\Omega^*$ of equal area; also note that,
since $u_*$ is radially decreasing, the ball
$\Omega^*$ is a superlevel set of $u_*$, as $E$ is of $u_\alpha$).

The meaning of \eqref{goal3} is now clear: given
any number $s>0$, 
the integral of a
function $u_\alpha$ of the kind
\eqref{not1}, over its superlevel set $E$ of measure $|E|=s$,
is \emph{maximized} when $\alpha=(0,1,\ldots,\newk-1)$, i.e.
when $u_\alpha=u_*$. We will prove this in a \emph{constructive way},
modifying the initial vector $\alpha=(k_1,\ldots,k_\newk)$ through a sequence
of moves: after each move (where
suitable entries of $\alpha$ are decreased by 1), the integral  of the new $u_\alpha$ 
(over its superlevel set of measure $s$) will increase, until
$\alpha$ is eventually turned into the optimal vector
$(0,1,\ldots,\newk-1)$.

\begin{proof}[Proof of \eqref{goal3}]
Consider an arbitrary vector $\alpha=(k_1,\ldots,k_\newk)$ with the $k_n$'s as in
\eqref{sort1}, and assume that $\alpha\neq (0,1,\ldots,\newk-1)$.
Due to \eqref{sort1}, denoting by $\newj$ the number of 
those subscripts
$n\in\{1,\ldots,\newk\}$ for which
$k_n=n-1$,  
we may either have $\newj \geq 1$,
  in which case
$k_n=n-1$ if $n\leq \newj$ and $k_n>n-1$ if $n>\newj$
(note that $\newj<\newk$ since $\alpha\neq (0,\ldots,\newk-1)$),
or we may have $\newj=0$,
in which case $k_n>n-1$ for every $n$.
Thus, 
the entries of the new vector
\begin{equation}
    \label{defnewa}
\newalpha:=
\begin{cases}
    (k_1-1,\ldots,k_\newk -1) & \text{if $\newj=0$}\\
    (0,1,\ldots,\newj-1,\,\,\, k_{\newj+1} -1,\ldots,k_\newk-1)
    & \text{if $\newj\geq 1$}
\end{cases}
\end{equation}
still form an increasing sequence, as did the
entries of $\alpha$ due to \eqref{sort1}.
\begin{lemma}\label{basiclemma}
    If $E_\alpha$ is any superlevel set of $u_\alpha$, and 
    $E_\newalpha$ is the superlevel set of
    $u_\newalpha$ such that
    $|E_\newalpha|=|E_\alpha|$, we have
\begin{equation}
\label{goal4}
\int_{E_\alpha} u_\alpha(z)\,dA(z)\leq  \int_{E_\newalpha} u_\newalpha(z)\,dA(z).
\end{equation}
Moreover, the inequality is strict if $\newj>0$, or if
$\newj=0$ and
$\newalpha=(0,1,\ldots,\newk-1)$.
\end{lemma}
\begin{proof}
Recalling \eqref{not1} and the
previous definition of $\newj$, we  can split
\begin{equation}
\label{split1}
    u_\alpha(z)=f(\pi |z|^2)+g(\pi|z|^2),\quad
f(\sigma):=e^{-\sigma}\sum_{n=0}^{\newj-1}  \frac{\sigma^{n}}{n!},
\quad
g(\sigma):=e^{-\sigma}\sum_{n=\newj+1}^\newk \frac{\sigma^{k_n}}{k_n!}
\end{equation}
(we agree that $f\equiv 0$ if $\newj=0$) and, similarly,
\begin{equation}
    \label{split2}
    u_\newalpha(z)=f(\pi |z|^2)+\newg(\pi|z|^2),\quad
\quad
\newg(\sigma):=e^{-\sigma}\sum_{n=\newj+1}^\newk \frac{\sigma^{k_n-1}}{(k_n-1)!}.
\end{equation}
   By Remark \ref{remsuplevel}, the structure of the set $E_\alpha$  
is   as in \eqref{structure}: therefore, letting $r_0=0$ to regard
the ball in \eqref{structure} as an annulus, 
   using polar coordinates and then letting 
   $\sigma=\pi r^2$
   we have
\begin{equation}
\label{polar}
\int_{E_\alpha} u_\alpha(z)\,dA(z)=
\sum_{i=0}^m
2\pi\int_{r_i}^{R_i}r 
\left(f(\pi r^2)+g(\pi r^2)\right)\,dr
=
\sum_{i=0}^m \int_{a_i}^{b_i}
\left( f(\sigma)+g(\sigma)\right)\,d\sigma,
\end{equation}
where we have set 
$a_i=\pi r_i^2$ and $b_i=\pi R_i^2$. Now we claim that
\begin{equation}
\label{bdrypoint}
g(a_i)-g(b_i)=f(b_i)-f(a_i)\leq 0\quad \forall i\in\{1,\ldots,m\}.
\end{equation}
This is obvious if $a_i=b_i$ (i.e. $r_i=R_i$ in \eqref{structure}, which corresponds to an empty annulus),
while if $a_i<b_i$ (i.e. $r_i<R_i$ and
the annulus is nonempty) the equality in \eqref{bdrypoint} follows from the 
splitting \eqref{split1} and fact that
$E_\alpha$ is a superlevel
set of $u_\alpha$, hence $u_\alpha$ is \emph{constant} along
$\partial E_\alpha$ (in particular, $u_\alpha$ is constant along the 
boundary of the $i$-th
annulus); finally, the inequality in \eqref{bdrypoint}
is \emph{strict} if $\newj>0$ because
$f(\sigma)$ is strictly decreasing in this case, 
whereas it is an equality if $\newj=0$ and $f(\sigma)\equiv 0$.
Combining \eqref{bdrypoint} with the fact that  $a_0=0$
and hence $g(b_0)-g(a_0)=g(b_0) \geq 0$ (again, 
equality holds when $\newj=0$, i.e. when $u_\alpha(0)=0$
and the ball in
\eqref{bdrypoint} is empty so that $R_0=b_0=0$), we finally have
\begin{equation}
    \label{ineq1}
\sum_{i=0}^m \bigl(g(a_i)-g(b_i)\bigr)
\leq 0,
\quad\text{with strict inequality if $\newj>0$.}
\end{equation}
Hence, observing that $g(\sigma)=\newg(\sigma)-g'(\sigma)$, 
integrating we have
\begin{equation}
\label{ineq2}
    \sum_{i=0}^m \int_{a_i}^{b_i}g(\sigma)\,d\sigma=
    \sum_{i=0}^m \int_{a_i}^{b_i}\newg(\sigma)\,d\sigma+
    \sum_{i=0}^m \bigl(g(a_i)-g(b_i)\bigr)
    \leq
    \sum_{i=0}^m \int_{a_i}^{b_i}\newg(\sigma)\,d\sigma,
\end{equation}
and the inequality is strict if $\newj>0$. Plugging this into
\eqref{polar} we get
\begin{equation}
    \label{ineq3}
\int_{E_\alpha} u_\alpha(z)\,dA(z)
\leq
\sum_{i=0}^m \int_{a_i}^{b_i}
\left( f(\sigma)+\newg(\sigma)\right)\,d\sigma
=\int_{E_\alpha} u_\newalpha(z)\,dA(z),
\end{equation}
where the last equality follows from \eqref{split2}, getting back
to polar coordinates with the change of variable $\sigma=\pi r^2$
(exactly as we got \eqref{polar} from \eqref{split1}, but in reverse order). 
From
the bathtub principle \eqref{bathtub} applied to $u_\newalpha$,
the last integral increases if we pass from $E_\alpha$
to the superlevel set $E_\newalpha$ of $u_\newalpha$ having the same measure, i.e.
\begin{equation}
\label{bathtub2}    
\int_{E_\alpha} u_\newalpha(z)\,dA(z)
\leq
\int_{E_\newalpha} u_\newalpha(z)\,dA(z),
\end{equation}
which combined with \eqref{ineq3} proves \eqref{goal4}.

Finally, assume that  \eqref{goal4} is an equality. Then equality
holds also in \eqref{ineq3} (hence $\newj=0$ otherwise
\eqref{ineq2} would be strict) and
in \eqref{bathtub2},
which implies
(as observed at the end of Remark \ref{remsuplevel}
relative to \eqref{bathtub})
that
$\left\vert E_\alpha \triangle E_\newalpha\right\vert=0$, i.e.,
that $E_\alpha$ is already
a superlevel set of $u_\newalpha$. To conclude the proof,
we now further assume that
$\newalpha=(0,1,\ldots,\newk-1)$ 
(hence $\alpha=(1,2,\ldots,\newk)$ by \eqref{defnewa} since $\newj=0$),
and show that this leads to a contradiction.
Indeed, with this assumption
$u_\newalpha(z)=u_*(z)$ (cf. \eqref{not1}) and hence,
as $u_*(z)$ is radially decreasing, $E_\newalpha$ is a ball.
Hence, from \eqref{not1} (with $k_n=n)$
and \eqref{eigenval_ball0} (with $\Omega=E_\alpha$)
we have
\[
\int_{E_\alpha} u_\alpha(z)\,dA(z)
=\sum_{n=1}^\newk \int_{E_\alpha} e^{-\pi|z|^2}|e_n(z)|^2\,dA(z)
=\lambda_2(E_\alpha)+\cdots+\lambda_{\newk+1}(E_\alpha),
\]
and,
similarly, 
\[
\int_{E_\alpha} u_\newalpha(z)\,dA(z)
=\sum_{n=0}^{\newk-1} \int_{E_\alpha} e^{-\pi|z|^2}|e_n(z)|^2\,dA(z)
=\lambda_1(E_\alpha)+\cdots+\lambda_{\newk}(E_\alpha).
\]
Since \eqref{ineq3} is now an equality, equating the last two sums
we find $\lambda_{\newk+1}(E_\alpha)=\lambda_1(E_\alpha)$, which
is impossible since, in the case of a ball, the eigenvalues in
\eqref{eigenval_ball0} are known to be \emph{simple} (see \cite{daubechies,seip}).
\end{proof}
With this lemma at hand, we can now prove \eqref{goal3}
(and hence also \eqref{eq:main estimate}, as explained after \eqref{goal1}).
Starting with
an arbitrary $\alpha\neq (0,1,\ldots,\newk-1)$ (otherwise
\eqref{goal3} is trivial), we can apply
\eqref{goal4} iteratively, replacing $\alpha$ with $\newalpha$ after
each iteration until eventually
$\newalpha=(0,1,\ldots,\newk-1)$, thus proving \eqref{goal3}
Now assume that equality occurs in \eqref{eq:main estimate}, hence also
in \eqref{goal3}. If $\alpha\neq (0,1,\ldots,\newk-1)$, then each
iteration of Lemma \ref{basiclemma} must yield equality in \eqref{goal4}
(otherwise in the end one would have strict inequality in \eqref{goal3}):
then, as claimed by the lemma, at each iteration
$\newj=0$ and $\newalpha\not=(0,1,\ldots,\newk-1)$. But this is
a contradiction, since the last iteration occurs precisely when
$\newalpha=(0,1,\ldots,\newk-1)$. Thus, equality in 
\eqref{goal3} occurs only if $\alpha=(0,1,\ldots,\newk-1)$,
i.e. when $u_\alpha=u_*$ and $E=\Omega_*$ is a ball. Tracing back to \eqref{goal1}
through \eqref{goal2}, equality in \eqref{goal1} is possible only 
if the initial $\Omega$ is a ball and (as discussed after \eqref{minmax1})
$F_1,\ldots,F_\newk$ span the same subspace as $e_0,\ldots,e_{\newk-1}$
(the first $\newk$ eigenfunctions of the ball $\Omega$). But since both
sets of functions are orthonormal in $\Fock$, this condition is equivalent to 
\eqref{unitU2}.
Hence Theorem \ref{th:main} is completely proved.
\end{proof}

\begin{remark}
It appears from the proof of Lemma \ref{basiclemma} that
\eqref{ineq3} reduces to an equality when $\newj=0$ (this
condition is equivalent to $u_\alpha(0)=0$, i.e. $k_1>0$
in \eqref{not1}, by \eqref{sort1}). Since the domain of integration
$E_\alpha$ is an arbitrary superlevel set of $u_\alpha$
(playing no apparent role for $u_\newalpha$),
equality in \eqref{ineq3} whenever $u_\alpha(0)=0$ is
somewhat surprising, even in the simplest case where $\newk=1$
in \eqref{not1} (i.e. $\alpha=(k_1)$ with $k_1>0$) and in \eqref{defnewa} 
(i.e. $\newalpha=(k_1-1)$). In this case, 
writing
$k$ in place of $k_1$ and recalling \eqref{monomials},
\eqref{ineq3} 
amounts to
\begin{equation}
\label{eqmon}
    \int_{E_k} e^{-\pi |z|^2} |e_k(z)|^2\,
    dA(z)
    =
    \int_{E_k} e^{-\pi |z|^2} |e_{k-1}(z)|^2\,
    dA(z)\quad k=1,2,\ldots
\end{equation}
where $E_k$ is an arbitrary superlevel set of the form
\begin{equation}
    \label{defEk}
    E_k=\left\{z\in\C\,|\,\, e^{-\pi |z|^2} |e_k(z)|^2 > t\right\}
    \quad\text{for some $t>0$.}
\end{equation}
When nonempty, the set $E_k$ is clearly an open annulus whose
inner and outer radii are the two positive solutions of the equation
\[
e^{-\pi r^2} \frac {\pi^k r^{2k}}{k!}=t
\]
(by \eqref{monomials}, two distinct solutions exist if and only if $t<e^{-k}k^k/k!$, otherwise
$E_k=\emptyset$). Equation \eqref{eqmon} (which one can easily  check 
integrating by parts) will play a central role in the next section.
\end{remark}

    \section{Consequences of the main result}\label{sec:consequences of the main result}
    
    Theorems \ref{teoPS} and \ref{th:main} have a number of interesting corollaries that we are going to state and prove in this section.

    The first corollaries are simple applications of 
     \emph{majorization theory} (see \cite{hardy_littlewood_polya}).
Given two vectors $x=(x_k)_{k=1}^n$ and $y=(y_k)_{k=1}^n$ in $\R^n$,
such that $x_k\geq x_{k+1}$ and $y_k\geq y_{k+1}$ ($1\leq k < n$), 
$x$ is said to be ``weakly majorized'' by $y$ (in symbols $x\prec_w y$)
if
\begin{equation}
\label{xsuccwy}
\sum_{k=1}^\newk x_k \leq \sum_{k=1}^\newk y_k,\quad
\forall \newk \in\{1,\ldots,n\}
\end{equation}
(``strong majorization'' or simply ``majorization'', 
in symbols $x\prec y$, occurs if, in addition, there is equality in \eqref{xsuccwy}
for $\newk=n$).
In these terms, the validity of the inequality in \eqref{ineqPS} 
for arbitrary $\newk\geq 1$ is equivalent to the weak majorizations
\begin{equation}
    \label{weakmaj}
    \left(\lambda_1(\Omega),\ldots,\lambda_n(\Omega)\right)
    \prec_w \left(\lambda_1(\Omega^*),\ldots,\lambda_n(\Omega^*)\right)
\quad\forall n\geq 1.
\end{equation}
By a well known theorem of Ky Fan \cite{kyfan},     
if $f:\R^n\to\R$ is a \emph{symmetric gauge function}
(i.e. $f$ is a norm on $\R^n$, $f(x_1,\ldots,x_n)=f(|x_1|,\ldots,|x_n|)$
and $f(x)=f(Px)$ for every $x\in\R^n$ and every permutation matrix $P$),
then $f(x)\leq f(y)$ whenever $x,y$ are vectors with positive entries
such that $x\prec_w y$. Thus, we obtain from Theorem \ref{teoPS}
and \eqref{weakmaj}:
\begin{corollary}
    \label{corgauge}
Let $\Omega\subset\C$ be a circularly symmetric
set of finite area $|\Omega|>0$, and let $\Omega^*$ be the ball
such that $|\Omega^*|=|\Omega|$. Then, for every $n\geq 1$ and every
symmetric gauge
function $f$ on $\R^n$, there holds
\begin{equation}
    \label{ineqgauge}
f \left(\lambda_1(\Omega),\ldots,\lambda_n(\Omega)\right)
    \leq f\left(\lambda_1(\Omega^*),\ldots,\lambda_n(\Omega^*)\right).
\end{equation}
\end{corollary}
A possible application is the case of    
    a weighted sum of the first $n$ eigenvalues:
    \begin{corollary}
Let $\Omega\subset\C$ be a circularly symmetric
set of finite area $|\Omega|>0$, and let $\Omega^*$ be the ball
such that $|\Omega^*|=|\Omega|$. If $n\geq 2$ and $t_1 \geq t_2 \geq \cdots \geq t_n > 0$,
then
        \begin{equation}\label{eq:inequality for weighted sum}
            \sum_{k=1}^n t_k\lambda_k(\Omega) \leq \sum_{k=1}^n t_k\lambda_k(\Omega^*),
        \end{equation}
and equality occurs if and only if $\Omega$ is equivalent
(up to a translation) to $\Omega^*$.
    \end{corollary}
The inequality in \eqref{eq:inequality for weighted sum} follows
easily from Corollary \ref{corgauge}, yet the following short proof
also enables us to characterize the cases of equality.
    \begin{proof}
Letting $\tau_{n}=t_n$ and $\tau_\newk=t_\newk-t_{\newk+1}$
for $1\leq \newk<n$,
there holds
\begin{equation}
    \label{eq33}
\sum_{k=1}^n t_k\lambda_k(\Omega)
=
\sum_{\newk=1}^n \tau_{\newk}\sum_{k=1}^\newk 
\lambda_k(\Omega)
\leq
\sum_{\newk=1}^n \tau_{\newk}\sum_{k=1}^\newk \lambda_k(\Omega^*)
=\sum_{k=1}^n t_k\lambda_k(\Omega^*)
\end{equation}
where we have used \eqref{ineqPS} with $\newk\in\{1,\ldots,n\}$ to estimate the inner sums (note $\tau_\newk\geq 0$). Since however $\tau_n>0$, if equality
holds in \eqref{eq33} then in particular equality must hold in \eqref{ineqPS} when
$\newk=n$, so $\Omega$ is equivalent to $\Omega^*$ by Theorem \ref{teoPS}.
    \end{proof}
\begin{remark}
The monotonicity assumption on the weights $t_k$ cannot be dropped. For
instance, if for some  $k_0<\newk$ we have $t_1 = \cdots = t_{k_0} < t_{k_0+1}$ then \eqref{eq:inequality for weighted sum} is false
for certain sets $\Omega$. Indeed, in this case one can check that the radial function
\[
u(z)=
\sum_{k=1}^\newk  t_k |e_k(z)|^2 e^{-\pi |z|^2}
\]
(whose integral on the ball $\Omega^*$, by \eqref{eigenval_ball0}, gives the right-hand side of
\eqref{eq:inequality for weighted sum}) is
radially \emph{increasing} in a neighborhood of $z=0$. Hence, if
$r>0$ is small enough and $\Omega$ is the annulus $2 r<|z|<\sqrt 5 \, r$
(so that $\Omega^*$ is the ball of radius $r$), then 
$\sup_{\Omega^*} u < \inf_{\Omega} u$ and hence
\[
\sum_{k=1}^n t_k\lambda_k(\Omega^*)
=
\int_{\Omega^*} u(z)\,dA(z) < 
\int_{\Omega} u(z)\,dA(z)\leq
\sum_{k=1}^n t_k\lambda_k(\Omega)
\]       
(the last inequality holds because the eigenfunctions of $T_\Omega$
are the monomials \eqref{monomials}).
\end{remark}
The following result is another immediate consequence of the weak majorizations \eqref{weakmaj}
(using \cite[Theorem 15.16]{simon_convexity}).
\begin{corollary}
Let $\Omega\subset\C$ be a circularly symmetric
set of finite area $|\Omega|>0$, and let $\Omega^*$ be the ball
such that $|\Omega^*|=|\Omega|$. Then, for every $n\geq 1$ and every
convex increasing function $\Phi:(0,1)\to\R$, there holds
\begin{equation}
    \label{ineqconvinc}
\sum_{k=1}^n \Phi(\lambda_k(\Omega))
\leq
\sum_{k=1}^n \Phi(\lambda_k(\Omega^*)).
\end{equation}
\end{corollary}
Letting $n\to\infty$ in \eqref{ineqconvinc} leads to inequalities of the 
kind \eqref{eq:inequality traces}, but under the extra assumption that $\Phi$
is increasing (the choice $\Phi(t)=t^p$ with $p\geq 1$).
In fact, the proof of Theorem \ref{teotrace} is more
delicate and requires the following extension of Karamata's
inequality:

\begin{proposition}\label{prop:Karamata inequality}
        Let $\{x_k\}, \{y_k\}\subset (0,1)$ be two decreasing sequences ($k\geq 1$) such that
        \begin{equation}
            \label{ass1}
            \sum_{k=1}^\newk x_k \leq \sum_{k=1}^\newk y_k\quad\forall\newk\geq 1,\quad\text{and}
            \quad\sum_{k=1}^{\infty} x_k = \sum_{k=1}^{\infty} y_k < \infty.
        \end{equation}
        Then, for every convex function $\Phi \colon (0,1) \to \R$ it holds
        \begin{equation}\label{eq:Karamata inequality}
            \sum_{k=1}^{\infty} \Phi(x_k) \leq \sum_{k=1}^{\infty} \Phi(y_k).
        \end{equation}
        Now assume, in addition, that $\Phi$ is  not an affine function on the interval
        $(0,y_1]$. If the series 
        $\sum \Phi(y_k)$ converges and equality occurs in \eqref{eq:Karamata inequality},  
        then the equality $\sum_{k=1}^\newk x_k = \sum_{k=1}^\newk y_k$ 
        occurs for at least one value of $\newk \geq 1$.
    \end{proposition}
The classical Karamata's inequality (see \cite[Theorem 108]{hardy_littlewood_polya})
states that $\sum_{k=1}^n\Phi(x_k)\leq \Phi(y_k)$, when $x=(x_k)$ and $y=(y_k)$
are vectors in $\R^n$ with decreasing entries such that $x\prec y$ (as defined after 
\eqref{xsuccwy}), and \eqref{eq:Karamata inequality} is a natural extension to the case 
of summable sequences (the conditions in \eqref{ass1} are the natural analogue
of the strong majorization $x\prec y$ for finite vectors). Now, while
the inequality 
\eqref{eq:Karamata inequality} is essentially already known (see \cite{filomat}), the statement
about the equality case is likely to be new (note that $\Phi$ is not assumed to be
strictly convex, but merely non-affine). Since it is not easily obtainable directly
from \eqref{eq:Karamata inequality}, a full proof of Proposition \ref{prop:Karamata inequality} is given in the Appendix.

\begin{proof}[Proof of Theorem \ref{teotrace}] Since $\Omega$ and $\Omega^*$ have
equal measure, 
due to \eqref{traccia} we have
\begin{equation}
    \label{equalsum}
\sum_{k=1}^\infty \lambda_k(\Omega)=
\sum_{k=1}^\infty \lambda_k(\Omega^*)<\infty
\end{equation}
and hence, since the inequality in \eqref{ineqPS} holds for every $\newk\geq 1$,
we can apply Proposition \ref{prop:Karamata inequality} to the two decreasing  sequences $x_k=\lambda_k(\Omega)$ and $y_k=\lambda_k(\Omega^*)$. Then \eqref{eq:inequality traces}
follows from \eqref{eq:Karamata inequality} and, in case of equality (with
finite sums and $\Phi$ not affine), by Proposition \ref{prop:Karamata inequality}
equality occurs in \eqref{ineqPS} for at least one values of $\newk\geq 1$,
and the last claim of  Theorem \ref{teoPS} implies
that $\Omega$ is a ball.
	\end{proof}
             \begin{remark}  
	Theorems \ref{teoPS} and \ref{th:main} can be rephrased in terms of the short-time Fourier transform with Gaussian window $\varphi(x) = 2^{1/4} e^{-\pi x^2}$, defined as
	\begin{equation*}
		\V f (x, \omega) = \int_{\R} f(t) \varphi(t-x) e^{-2\pi i \omega t} \, dt, \quad (x,\omega) \in \R^2,\ f \in L^2(\R).
	\end{equation*}
    Indeed, it is well known (cf. \cite[Section 3.4]{grochenig}) that the Bargmann transform
    \begin{equation*}
        \mathcal{B}f(z)\coloneqq 2^{1/4} \int_{\R} f(t) e^{\pi i t z - \pi t^2 - \pi |z|^2/2}\, dt, \quad z = x+i\omega \in \C
    \end{equation*}
    is a unitary operator from $L^2(\R)$ on $\Fock$ and is related to the short-time Fourier transform via the relation
    \begin{equation*}
        \V f(x,-\omega) = e^{\pi i x \omega} e^{- \pi |z|^2/2} \mathcal{B}f(z), \quad z = x + i \omega \in \C.
    \end{equation*}
   Moreover, the Hermite functions
    \begin{equation*}
            h_k(t) = \frac{2^{1/4}}{\sqrt{k!}}\left(-\frac{1}{2\sqrt{\pi}}\right)^k e^{\pi t^2} \frac{\mathrm{d}^k}{\mathrm{d}t^k}(e^{-2\pi t^2}), \quad t \in \R,\ k \in \N
    \end{equation*}
    are mapped by $\mathcal{B}$ into the normalized monomials in $\Fock$. 
    
    We can then consider the so-called time-frequency localization operator, or anti-Wick operators (see \cite{daubechies}), defined as \[
    L_{\Omega} \coloneqq \V^* \chi_{\Omega} \V \colon L^2(\R) \to L^2(\R),
    \]
    where $\Omega \subset \R^2$ is a set of finite area. With the natural identification $\R^2 \simeq \C$ it is easy to see that, if $\Omega$ is radial (i.e., circularly symmetric around the origin), then $L_{\Omega} = \mathcal{B}^* T_{\Omega} \mathcal{B}$ and therefore $L_{\Omega}$ has the same eigenvalues of $T_{\Omega}$. Indeed, if $F=\mathcal{B}f$, $G=\mathcal{B}g$, with $f,g\in L^2(\R)$, we have 
    \begin{align*}
        \langle T_\Omega F,G\rangle_{\Fock}&=\int_\Omega F(z) \overline{G(z)} e^{-\pi|z|^2} \,dA(z)=\int_\Omega \mathcal{B}f(z)\overline{\mathcal{B}g(z)} e^{-\pi|z|^2}\, dA(z)\\
        &=\int_\Omega \mathcal{V}f(x,-\omega)\overline{\mathcal{V}g(x,-\omega)}\, dx\,d\omega=\langle L_\Omega f,g\rangle_{L^2}.
    \end{align*}
    With this in mind, Theorems \ref{teoPS} and \ref{th:main} can be easily rephrased in terms of the short-time Fourier transform and the operator $L_\Omega$. 
\end{remark}
    \section{Maximizing the second eigenvalue}\label{sec:maximizing the second eigenvalue}
    
In this section we address the problem of maximizing the second eigenvaule $\lambda_2(\Omega)$ of $T_\Omega$, among circularly symmetric sets $\Omega$
of given area. In particular, we will prove Theorem \ref{teolambda2} and Proposition \ref{prononrad}.
\begin{proof}[Proof of Theorem \ref{teolambda2}]
Let $\Omega$ be as in the statement of the theorem and let $s=|\Omega|$
denote its area. By circular symmetry (which by Remark \ref{reminv}
we may assume to be \emph{around the origin}), the eigenfunctions
of $T_\Omega$ are the monomials in \eqref{monomials}, and 
the eigenvalues
$\lambda_n(\Omega)$ are given (in some unknown order but with the correct
multiplicities)
by the integrals in \eqref{eq:expression eigenvalues}. 
In particular, letting for simplicity
\[
u_k(z):=e^{-\pi |z|^2}|e_k(z)|^2=e^{-\pi|z|^2}\frac{\pi^k |z|^{2k}}{k!},
\quad k=0,1,2,\ldots
\]
we have that
\begin{equation}
    \label{charl2}
\lambda_2(\Omega)\quad\text{is the second largest value of
$\,\,\,\int_\Omega u_k(z)\,dA(z),\,\,$ as $k\geq 0.$}
\end{equation}
Denoting by $E_k^s$ the superlevel set of $u_k$ of measure $s$,
that is, a set as in \eqref{defEk} where $t$ is now chosen in such a way
that $|E_k^s|=s$ (cf. the last sentence in Remark \ref{remsuplevel}),
we have from \eqref{eqmon}
\begin{equation}
    \label{iter1}
    \int_{E_k^s} u_k(z)\,dA(z)=
    \int_{E_k^s} u_{k-1}(z)\,dA(z)
    <\int_{E_{k-1}^s} u_{k-1}(z)\,dA(z)\quad\forall k\geq 1,
\end{equation}
where the inequality follows from the bathtub principle \eqref{bathtub}
and is strict because it holds $|E_k^s\triangle E_{k-1}^s|>0$
(note that $E_0^s$ is the ball of area $s$, while the other
$E_k^s$ are pairwise distinct annuli, as described after \eqref{defEk}).

If, starting from some $k\geq 2$, we iterate \eqref{iter1}
$k\,$--$1$ times, we obtain the inequality
\begin{equation}
\label{iter2}
\int_{E_k^s} u_k(z)\,dA(z)< \Lambda:=
    \int_{E_1^s} u_{1}(z)\,dA(z)
    =\int_{E_{1}^s} u_{0}(z)\,dA(z)\quad\forall k\geq 2
\end{equation}
(the former equality serves to define $\Lambda$, whereas the latter is a further application of \eqref{eqmon} with $k=1$),
and combining the bathtub principle \eqref{bathtub} with \eqref{iter2}
we find
\begin{equation}
\label{stima1}
    \int_\Omega u_k(z)\,dA(z)\leq \int_{E_k^s} u_k(z)\,dA(z)<
    \Lambda\quad\forall k\geq 2.
\end{equation}
Moreover, now using only the bathtub principle, we have
\begin{equation}
\label{stima2}
    \int_\Omega u_1(z)\,dA(z)\leq \int_{E_1^s} u_1(z)\,dA(z)=
    \Lambda\quad\text{(with strict inequality if $|\Omega\triangle E_1^s|>0$).}
\end{equation}
Now, a combination of \eqref{charl2} with \eqref{stima1} and \eqref{stima2}
proves the following: $\lambda_2(\Omega)\leq\Lambda$, and the inequality is
strict unless $|\Omega\triangle E_1^s|=0$ (i.e., unless $\Omega=E_1^s$).
When $\Omega=E_1^s$, this statement reveals that $\lambda_2(E_1^s)=\Lambda$
and, since the only integral not estimated in \eqref{stima1} or \eqref{stima2}
is when $k=0$, 
by the equality of the two integrals in \eqref{iter2}
we see that also $\lambda_1(E_1^s)=\Lambda$. Summing up: the annulus $E_1^s$
maximizes $\lambda_2(\Omega)$, among all
circularly symmetric sets $\Omega$ of measure $s$. Since clearly
$E_1^s=A_s$ as defined in \eqref{defAs}, to complete the proof
of Theorem \ref{teolambda2} it remains to prove that
the inner and outer radii of $A_s$ (as defined in \eqref{defAs})
are those given by \eqref{defradii}, and to check 
the explicit value of $\lambda_2(A_s)$ in \eqref{inequality lambda 2}.

The set $A_s$ in \eqref{defAs} is clearly an annulus
$\{r<|z|<R\}$ whose area $\pi(R^2-r^2)$ must equal $s$, and 
we must also have $u_1(r)=u_1(R)=t$. 
Thus, since $u_1(\rho)=e^{-\pi\rho^2}\pi\rho^2$,
setting $a=\pi r^2$ and $b=\pi R^2$ the annulus $A_s=E_1^s$
is determined by the system two equations $b-a=s$ and $a e^{-a}=b e^{-b}$.
Since the only (positive) solution is given by
$a:=s/(e^s-1)$ and $b:=s e^s/(e^s-1)$, \eqref{defradii} is proved.
Finally, since $u_0(z)=e^{-\pi |z|^2}$ and $A_s=E_1^s$,
we have from \eqref{iter2} that the explicit value of $\lambda_2(A_s)=\Lambda$
is given by
\begin{align*}
    \int_{A_s} e^{-\pi |z|^2}\,dA(z)
    =2\pi\int_r^R \rho e^{-\pi\rho^2}\,d\rho=
    \int_a^b e^{-\sigma}\,d\sigma=e^{-a}\left(1-e^{a-b}\right)
    =e^{-a}(1-e^{-s})
\end{align*}
and, since $a=s/(e^s-1)$, the equality
in \eqref{inequality lambda 2} is proved.
\end{proof}

\begin{remark}
          There is no analogue of Theorem \ref{teolambda2} for the third eigenvalue, that is, it is not true that the annulus
          $E_2^s$ (superlevel set of $u_2$ of measure $s$) maximizes $\lambda_3(\Omega)$, among
          circularly symmetric sets of area $s$. Indeed, since
          (generalizing \eqref{charl2})
          $\lambda_n(E_2^s)$ is the \emph{$n$-th largest} value among
          the integrals $\int_{E_2^s} u_k(z)\,dA(z)$ as $k\geq 0$, it
          is easy to check that
\[
\lambda_1(E_2^s)=\int_{E_2^s} u_2(z)\,dA(z)
=\int_{E_2^s} u_1(z)\,dA(z)=
\lambda_2(E_2^s)>\lambda_3(E_2^s)=\int_{E_2^s} u_0(z)\,dA(z)
\]          
(the other integrals with $k\geq 3$ provide smaller values).
Now, given  $\eps\in(0,s)$, consider
the competitor $\Omega_\eps:=A_\eps\cup B_\eps$, 
where $A_\eps\subset E_2^s$
is any  smaller annulus
of area $s-\eps$, and $B_\eps$ is the ball of area $\eps$
(we choose $\eps$ so small that $A_\eps\cap B_\eps=0$, in such a way that
$|\Omega_\eps|=s$).
Since $e_0(z)=e^{-\pi |z|^2}$ is radially decreasing,
replacing $E_2^s$ with $\Omega_\eps$ we will have
\[
\int_{\Omega_\eps} u_0(z)\,dA(z)>
\int_{E_2^s} u_0(z)\,dA(z)=\lambda_3(E_2^s)
\]
(the increase being of order $O(\eps)$)  whereas
the two integrals of
$u_1$ and $u_2$ will decrease a bit, but (for small enough $\eps$)
each of them will remain larger than $\lambda_3(E_2^s)$.
In other words, we will have
\[
\int_{\Omega_\eps} u_k(z)\,dA(z)>\lambda_3(E_2^s),
\quad k\in\{0,1,2\},
\]
and since each of these three integrals
is an eigenvalue of $T_{\Omega_\eps}$, we will
have $\lambda_3(\Omega_\eps)>\lambda_3(E_2^s)$.

Finally,  a refinement of this argument shows that, for $k\geq 3$,
the annulus $E_{k-1}^s$ (superlevel set of $u_{k-1}$ of measure
$s$) is not a maximizer of $\lambda_k(\Omega)$.
\end{remark}

We now recall some well known facts that
will be used in the proof of Proposition \ref{prononrad}
(see \cite{zhu_book} for
more details).
As  quickly mentioned after \eqref{concineq1},
the functions
$K_w(z)=e^{\pi\overline w z}$ are the reproducing kernels
of the Fock space $\Fock$, that is,
\begin{equation}
    \label{repk}
    F(w)=\langle F,K_w\rangle_\Fock =
    \int_{\C} F(z) e^{\pi w \overline{z}} e^{-\pi |z|^2}\,dA(z),
    \quad\forall F\in\Fock,\quad\forall w\in\C
\end{equation}
from which,
choosing $F=K_w$ and then $F=K_{-w}$, one can compute the quantities
\begin{equation}
    \label{normKw}
\Vert K_w\Vert_\Fock^2=e^{\pi |w|^2},\qquad
\langle K_{-w},K_w\rangle_\Fock=e^{-\pi |w|^2}.
\end{equation}
Finally, since equality is achieved in \eqref{concineq1} 
(as mentioned thereafter) when $F=K_w/||K_w||_\Fock$ and $\Omega$ is any ball centered at $w$, we have using \eqref{normKw}
\begin{equation}
    \label{Kwball}
\int_B e^{-\pi |z|^2} \vert K_w(z)\vert^2\,dA(z) = e^{\pi |w|^2}\bigl(
1-e^{-|B|}\bigr)
\quad\text{for every $B$ ball centered at $w$.}
\end{equation}

\begin{proof}[Proof of Proposition \ref{prononrad}]
Given a radius $r>0$, consider the domain $\Omega=B^+\cup B^-$, 
where $B^\pm$ denote the two balls of radius $r$ centered at $\pm w$,
and $w\in\R$ satisfies
$w>r$ so that $B^+\cap B^-=\emptyset$. The idea of the proof is that, 
when $w\gg r$,
the localization operator $T_\Omega$ (acting on $\Fock$) essentially behaves as a direct sum
$T_{B^+}\oplus T_{B^-}$ acting on $\Fock\oplus \Fock$ and hence, when $w\to\infty$, its spectrum tends to split as
$\sigma(T_\Omega)\sim \sigma(T_{B^+})\cup \sigma(T_{B^-})$
so that, in the limit, $\lambda_1(\Omega)=\lambda_2(\Omega)=\lambda_1(B^\pm)$. 
In other words, since letting $s=|\Omega|$ (i.e. $s=2\pi r^2$)
one has $|B^\pm|=s/2$ and 
hence $\lambda_1(B^\pm)=1-e^{-s/2}$ by \eqref{eigenval_ball0},
one expects that
\begin{equation}
    \label{as2}
\lambda_2(\Omega)\sim 1-e^{-s/2}\quad\text{as $w\to\infty$.}
\end{equation}
This asymptotics would complete the proof:
since by a Taylor expansion there holds
\begin{equation}
    \label{as3}
1-e^{-s/2} > e^{s/(1-e^s)}\left(1-e^{-s}\right) \quad\text{for small enough $s>0$,}
\end{equation}
we see that \eqref{tesiprop} would follow from \eqref{as2}, provided
$s$ is small enough and $w$ very large.

In fact, instead of \eqref{as2},
 we will 
prove
the one-sided  bound
\begin{equation}
    \label{lowerb}
\liminf_{w\to\infty}\,\,\lambda_2(\Omega)\geq 1-e^{-s/2}
\quad\text{for every fixed radius $r$,}
\end{equation}
which combined with \eqref{as3} proves \eqref{tesiprop}:
it will suffice to first fix the radius $r$ (hence also
the measure $s=2\pi r^2$)
small enough, so that
the inequality in \eqref{as3} holds true, and then choose
$w$ large enough,  to guarantee that the resulting
$\Omega$ satisfies
\eqref{tesiprop}.

So let us prove $\eqref{lowerb}$, with $\Omega=B^+\cup B^-$ and
$B^\pm$ balls of radius $r$ centered at $\pm w$, as
described above.
From the max-min characterization of the eigenvalues, we have
    \begin{equation}
        \label{maxmin} \lambda_2(\Omega)=\max_{S_2}\min_{F\in S_2\setminus\{0\}}\frac 
        {\int_\Omega e^{-\pi |z|^2} |F(z)|^2\, dA(z)}{\Vert F\Vert_{\Fock}^2},
    \end{equation}
    where the max is over all subspaces $S_2\subset\Fock$ of dimension two. Thus,
considering in particular the subspace $S_2$ spanned by the two reproducing kernels $k^\pm(z):= e^{\pm \pi wz}$, 
we have from \eqref{maxmin}
    \begin{equation}
        \label{maxmin2} \lambda_2(\Omega)\geq \min_{
        \begin{smallmatrix}
        a,b\in\C\\  (a,b)\not=(0,0)\end{smallmatrix}}\frac 
        {\int_\Omega e^{-\pi |z|^2} \left| a k^+(z)+b k^-(z)\right|^2\, dA(z)}{\Vert a k^+ +b k^-\Vert_{\Fock}^2}.
    \end{equation}
Using \eqref{normKw} (note $k^+=K_w$ and $k^-=K_{-w}$), one can easily compute
\begin{equation}
    \label{estden}
\Vert a k^+ +b k^-\Vert_{\Fock}^2=\left(|a|^2 +|b|^2\right) e^{\pi w^2}+2 e^{-\pi w^2} \Re a \overline  b
\leq
\left(|a|^2 +|b|^2\right) \left(e^{\pi w^2}+ e^{-\pi w^2}\right).
\end{equation}
On the other hand, writing $k^\pm$ and $d\nu$ as abbreviations 
for $k^\pm(z)$ and  
$e^{-\pi|z|^2}\,dA(z)$ in the integral,
since by \eqref{Kwball} $\int_{B^\pm}|k^\pm|^2\,d\nu=e^{\pi w^2}(1-e^{-s/2})$, we have 
\begin{equation*}
\begin{aligned}
    \int_\Omega  \left| a k^++b k^-\right|^2 d\nu
  &  > |a|^2\int_{B^+}   \left| k^+\right|^2 d\nu+
|b|^2\int_{B^-}   \left| k^-\right|^2d\nu -2 \Re  a\overline{b}\int_{\Omega} 
 k^+\,\overline {k^- }\,d\nu\\
& \geq
\left(|a|^2 +|b|^2\right)\left(e^{\pi w^2}\left(1-e^{-s/2}\right)
-\int_{\Omega} 
\left\vert k^+\,\overline {k^- }\right\vert\,d\nu\right)
\end{aligned}
\end{equation*}
and, since when $z\in\Omega$ we have $| k^+(z)\overline {k^-(z) }|=|e^{\pi w(z-\overline z)}|=|  e^{2\pi w i \Im z}|=1$, there holds
\begin{align*}
    \int_{\Omega}\left|   k^+\overline {k^- }\right|\,d\nu
     &=   \int_{\Omega}  \,d\nu
    =2 \int_{B^+} e^{-\pi (x^2+y^2)}\,dxdy\\ &\leq 2  |B^+|\, e^{-\pi(w-r)^2} 
    < 2\pi r^2 e^{-\pi w^2+2\pi wr }.
\end{align*}
Plugging this into the previous estimate, we finally obtain
\begin{align}
\label{estnum}
    \int_\Omega  \left| a k^+(z)+b k^-(z)\right|^2 \, d\nu
    >
    \left(|a|^2 +|b|^2\right)\left(e^{\pi w^2}\left(1-e^{-s/2}\right)
-2\pi r^2 e^{-\pi w^2+2\pi wr }
\right).
\end{align}
It is now easy to see that \eqref{estnum} and \eqref{estden},  
combined with
\eqref{maxmin2}, yield \eqref{lowerb}.
\end{proof}

\section{Extension to a general symbol}\label{sec:extension to a general symbol}
	The class of Toeplitz operators, considered in the previous sections, can be easily extended to include symbols that are not necessarily characteristic functions. Namely, consider a measurable function $\sigma \colon \C \rightarrow \C$, say $\sigma\in L^p(\C)=L^p(\C,dA)$, for some $p\in [1,+\infty]$, where, as usual, $dA$ denotes the Lebesgue measure in $\C\simeq\R^2$. Let $P$ be the projection operator from $L^2(\C, e^{-\pi |z|^2} dA(z))$ onto $\Fock$. We then define the Toeplitz operator \[
    T_{\sigma} := P \sigma \colon \Fock \rightarrow \Fock.
    \]
    When $\sigma=\chi_\Omega$, for some measurable subset $\Omega\subset\C$,  we recapture the operator $T_\Omega$ considered so far.
    The function $\sigma$ is usually referred to as \emph{symbol}, or \emph{weight} of the operator. It is well known that if $\sigma \in L^p(\C)$ for some $p \in [1,+\infty)$ then $T_\sigma$ is bounded and compact in $\Fock$ (\cite{daubechies,wong}). Moreover if $\sigma$ is, in addition, real-valued and nonnegative, then $T_\sigma$ is also self-adjoint and nonnegative, and therefore its spectrum consists of a sequence of eigenvalues $\lambda_1(\sigma) \geq \lambda_2(\sigma) \geq \cdots \geq 0$, converging to zero. In this section, we are going to study the problem of maximizing the sum $\sum_{k=1}^{\newk} \lambda_k(\sigma)$ when $\sigma$ is radial. 

    Observe that, in terms of the associated quadratic form, if $F\in \Fock$, we have
\[
\langle T_\sigma F,F\rangle_\Fock=\int_{\C}\sigma(z)|F(z)|^2 e^{-\pi|z|^2}\, dA(z).
\]

For a nonnegative symbol $\sigma \in L^p(\C)$ ($p<\infty$),
by the so called ``layer cake'' representation
\[
\sigma(z)=\int_0^\infty \chi_{\{\sigma>t\}}(z)\,dt,\quad z\in\C,
\]
we obtain, for an arbitrary $F\in\Fock$, that
\begin{equation}
\begin{aligned}
    \label{repF}
&\int_\C \sigma(z) |F(z)|^2 e^{-\pi |z|^2}\,dA(z)
=\int_\C \left(\int_0^\infty \chi_{\{\sigma>t\}}(z)\,dt
\right)|F(z)|^2 e^{-\pi |z|^2}\,dA(z)\\
&=\int_0^\infty \left(\int_{\{\sigma>t\}} |F(z)|^2 e^{-\pi |z|^2}\,dA(z)
\right)\,dt.
\end{aligned}
\end{equation}
Resuming the short notation $T_\Omega$ for $T_{\chi_\Omega}$ used
in the Introduction, \eqref{repF} becomes
\begin{equation}
\label{repF2}
\langle T_\sigma F,F\rangle 
=\int_0^\infty \langle T_{\{\sigma>t\}} F,F\rangle\,dt,
\end{equation}
which reveals $T_\sigma$ as a \emph{mixture} of the localization
operators $T_{\{\sigma>t\}}$ on the superlevel sets of $\sigma$. When $\sigma$ is radial,
this allows us to prove that 
the sum of the first $\newk$ eigenvalues increases, if we replace $\sigma$ with its decreasing rearrangement $\sigma^*$ defined by
    \begin{equation}
    \label{defFstar}
        \sigma^*(z) := \int_0^{\infty} \chi_{\{\sigma>t\}^*}(z) \, dt,\quad
        z\in\C,
    \end{equation}
    where $\{\sigma>t\}^*$ is the open ball of center 0 and the
    same measure as $\{\sigma>t\}$
(it is well known that $\sigma^*$ is radially decreasing and equimeasurable with $\sigma$, in particular $\|\sigma^*\|_p = \|\sigma\|_p$, see e.g. \cite{hardy_littlewood_polya}).

Defining for $\newk\geq 1$ the functions
	\begin{equation}
    \label{defGk}
		G_\newk(s) = \newk - \sum_{k=0}^{\newk-1} (\newk-k)\dfrac{s^k}{k!}e^{-s}
	\end{equation}
(these are functions appearing in the right-hand side of \eqref{eq:main estimate}), we have the following result. 

\begin{theorem}\label{teoF1}
		Assume $\sigma \in L^p(\C)$ ($1\leq p<\infty$) is a radial nonnegative symbol 
        and let $\sigma^*$ be its
        decreasing rearrangement as in \eqref{defFstar}. Then, for every $\newk\geq 1$
		\begin{equation}\label{rearr}
			\sum_{k=1}^{\newk} \lambda_k(\sigma) \leq \sum_{k=1}^{\newk} \lambda_k(\sigma^*) = \int_0^{\infty} G_\newk(\mu(t)) \, dt,
		\end{equation}
		where $\mu(t) = |\{\sigma>t\}|$ is the distribution function of $\sigma$. Moreover, equality in \eqref{rearr} is achieved if and only if $\sigma(z)=\sigma^*(z)$ 
        for a.e.  $z\in\C$. Finally, every eigenvalue $\lambda_n(\sigma^*)$
        ($n\geq 1$) is simple, and the corresponding eigenfunction 
        is the monomial $e_{n-1}(z)$ in \eqref{monomials}.
\end{theorem}
Here, as usual, $\sigma$ being \emph{radial} means that $\sigma(z)=\sigma(|z|)$.
 By \eqref{repF}, however, Remark \ref{reminv} carries over from
sets $\Omega\subset\C$ to arbitrary symbols $\sigma\in L^p$, hence the theorem
remains
valid if $\sigma$ is radial around some center $z_0\in\C$
(with equality when $\sigma(z)=\sigma^*(z-z_0)$).

\begin{proof} 
		Since $\sigma$ is radial, it is known (see \cite{daubechies}) that  the eigenfunctions of $T_\sigma$ are still the monomials 
        \eqref{monomials}
         and its eigenvalues are therefore given (in some order) by
        the numbers
        \begin{equation}
\label{eigenvaluesF}
\langle T_\sigma \,e_k,e_k\rangle_\Fock=        \int_{\C} 
\sigma(z)\,|e_k(z)|^2 e^{-\pi|z|^2} \, dA(z), \quad k=0,1,2,\ldots
\end{equation}
which, in particular, is true for $T_{\sigma^*}$. 
In this case,
however, using
\eqref{repF} (written for $\sigma^*$ first with $F=e_{k}$, $k\geq 1$,
and then with $F=e_{k-1}$) we can compare
the Rayleigh quotients
\begin{align*}
    &\langle T_{\sigma^*} \,e_{k},e_{k}\rangle_\Fock
    =
    \int_0^\infty \langle T_{\{\sigma^*>t\}} e_{k},e_{k}\rangle\,
    dt
    =
    \int_0^\infty \lambda_{k+1}(\{\sigma^*>t\}) \,    dt\\
&\qquad  >
\int_0^\infty \lambda_{k}(\{\sigma^*>t\}) \,    dt=
    \int_0^\infty \langle T_{\{\sigma^*>t\}} e_{k-1},e_{k-1}\rangle\,dt
    = \langle T_{\sigma^*} \,e_{k-1},e_{k-1}\rangle_{\Fock},
\end{align*}
having used \eqref{eigenval_ball0} with $\Omega=\{\sigma^*>t\}$
(which is a ball) and the fact that, for balls, every eigenvalue
is simple (see \cite{daubechies,seip}). This monotonicity reveals that, similarly
to \eqref{eigenval_ball0},
\begin{equation}
    \label{repr2}
\lambda_k(T_{\sigma^*})=\langle T_{\sigma^*} \,e_{k-1},e_{k-1}\rangle_\Fock
=\int_0^\infty \lambda_{k}(\{\sigma^*>t\}) \,    dt,
\end{equation}
thus proving the last claim of
the theorem. In addition, since $\{\sigma^*>t\}$ is a ball and by equimeasurability $\vert 
\{\sigma^*>t\}\vert=\mu(t)$, the last equality in \eqref{rearr} follows
from \eqref{repr2}, \eqref{ineqPS} (last equality, with $\Omega=\{\sigma>t\}$)
and \eqref{defGk}.

Then, getting back to $\sigma$ and the relative eigenvalues \eqref{eigenvaluesF},
let $k\to i_k$ be a bijection such that $e_{i_k}$ is an eigenfunction of $T_\sigma$ relative to the $k$-th
eigenvalue $\lambda_k(\sigma)$.
Letting $F=e_{i_k}$ in
\eqref{repF2} and summing,  we have
\begin{equation}
\label{ali1}
\begin{aligned}
&\sum_{k=1}^{\newk} \lambda_k(\sigma) =
\sum_{k=1}^{\newk} 
\langle T_\sigma e_{i_k},e_{i_k}\rangle
=
\int_0^{\infty} 
\sum_{k=1}^{\newk} 
\langle T_{\{\sigma>t\}} e_{i_k},e_{i_k}\rangle\,dt
\leq
\int_0^{\infty} 
\sum_{k=1}^{\newk} 
\lambda_k(\{\sigma>t\}) \,dt
\\
&\leq
\int_0^{\infty} 
\sum_{k=1}^{\newk} 
\lambda_k(\{\sigma>t\}^*) \,dt
=\int_0^{\infty} 
\sum_{k=1}^{\newk} 
\lambda_k(\{\sigma^*>t\}) \,dt=
\sum_{k=1}^{\newk} 
\lambda_k(T_{\sigma^*})
\end{aligned}
\end{equation}
where the two inequalities follow, respectively, from \eqref{minmax1}
and from \eqref{ineqPS}, and the last inequality follows from \eqref{repr2}.
This proves the inequality in \eqref{rearr}, equality being
possible only if equality holds, for a.e. $t$, also in \eqref{ineqPS}
when $\Omega=\{\sigma>t\}$. As claimed after \eqref{ineqPS},
this implies that $\{\sigma>t\}$ is  a ball 
for a.e. $t$ (hence for \emph{every} $t$, 
as $\{\sigma>t\}=\bigcup_{\tau>t} \{\sigma>t\}$), that is, $\sigma=\sigma^*$ almost everywhere.
Conversely, when $\sigma=\sigma^*$, the sets $\{\sigma>t\}$ are balls (hence there is
equality in \eqref{minmax1} and in \eqref{ineqPS} with
$\Omega=\{\sigma>t\}$) and by \eqref{repr2} the bijection $k\to i_k$
is simply $i_k=k-1$. All in all, when $\sigma=\sigma^*$ there is equality in \eqref{ali1} and hence also in \eqref{rearr}.
\end{proof}

The integral in \eqref{rearr} provides an optimal upper bound,
in terms of the distribution function $\mu(t)=|\{\sigma>t\}|$,
for the sum of the first $\newk$ eigenvalues of $T_\sigma$. It is then natural to investigate those symbols $\sigma\in L^p(\C)$ that maximize this sum,
under an $L^p$-norm constraint on $\sigma$. When $\newk=1$ (maximization
of $\lambda_1(\sigma)$) this problem was studied in  \cite{nicolatilli_norm, riccardi_optimal_estimate}, allowing for several  simultaneous norm constraints. Here, when $\newk>1$, combining Theorem \ref{teoF1}
with the techniques introduced in \cite{nicolatilli_norm}, we prove
the following result (where, for simplicity, we restrict ourselves to the case  of one $L^p$ constraint with $p>1$).

	\begin{theorem}\label{teoF2}
		Let $p > 1$, $B>0$ and $\newk \geq 1$. Then, for all radial nonnegative symbols $\sigma \in L^p(\C)$ such that $\|\sigma\|_p \leq B$, it holds
		\begin{equation}\label{rearr2}
			\sum_{k=1}^{\newk} \lambda_k(\sigma) \leq \sum_{k=1}^{\newk} \lambda_k(\sigma_p),
		\end{equation}
        where the optimal symbol $\sigma_p \in L^p(\C)$ is given by
            \begin{equation}\label{defFp}
                \sigma_p(z) = \frac B {c_p} \left( \sum_{k=0}^{\newk-1} \frac{(\pi |z|^2)^k}{k!}\right)^{\frac{1}{p-1}}e^{-\frac{\pi |z|^2}{p-1}}, \quad z \in \C,
            \end{equation}
        where $c_p>0$ is defined by 
        $c_p^p = \int_0^{\infty} (G_\newk'(s))^{\frac{p}{p-1}}\,ds$, and $G_\newk$ is the function in \eqref{defGk}. Moreover, equality in \eqref{rearr2} is achieved if and only if $\sigma(z)=\sigma_p(z)$ for a.e. $z \in \C$.
        \end{theorem}
\begin{proof}
If $\sigma$ is nonnegative and $\Vert \sigma\Vert_{L^p}\leq B$,
then its distribution function $\mu(t) = |\{\sigma> t\}|$ clearly
belongs to the class of functions
        \begin{equation*}
            \mathcal{C} := \left\{ \nu \colon (0,+\infty) \to [0,+\infty)\quad
            \text{such that $\,\,\nu$ is decreasing and  $\,\,p\int_0^{\infty} t^{p-1} \nu(t) \, dt \leq B^p$} \right\}.
        \end{equation*}
Therefore, 
we have from \eqref{rearr}
\begin{equation}
\label{varb}
            \sum_{k=1}^{\newk} \lambda_k(\sigma) \leq \int_0^{\infty} G_\newk(\mu(t)) \, dt
            \leq
                        \sup_{\nu \in \mathcal{C}} \int_0^{\infty} G_\newk(\nu(t)) \, dt.
\end{equation}
The idea of the proof is to show that the supremum in \eqref{varb}
is  attained by a unique $\nu\in {\mathcal C}$ which 
is the distribution function of the admissible symbol $\sigma_p$ 
defined in \eqref{defFp}.
When $\newk=1$ this was accomplished in
\cite{nicolatilli_norm}, under the additional
constraint that $\Vert\sigma\Vert_{L^\infty}\leq A\in (0,\infty]$
(letting $A=\infty$ makes this constraint inactive, as in our case).
Following step by step the proof of
\cite[Theorem 3.4]{nicolatilli_norm}, however,
reveals that the very same argument (with the same $B$
and in the particular case 
where
$A=\infty$,
but with the function $G=G_1$ replaced by $G_\newk$)
can 
be used here (the only relevant properties
of $G=G_1$ used in \cite{nicolatilli_norm} are smoothness, 
strict monotonicity and concavity, which $G_\newk$ satisfies).
In this way,  
one proves that the supremum in  \eqref{varb} 
is attained by a unique
$\mu_p\in {\mathcal C}$, that is implicitly defined by the condition that
        \begin{equation}\label{defmup}
            \frac B {c_p}G_\newk'(\mu_p(t))^{\frac 1{p-1}} = t\quad\text{if} \quad 0<t < \frac B {c_p},
        \end{equation}
whereas $\mu_p(t)=0$ for $t \geq B/c_p$. 
It is easy to see that $\mu_p$ is the
distribution function of $\sigma_p$. Indeed, 
as $\sigma_p$ is (strictly) radially decreasing, taking any $z\in\C$
and letting $t=\sigma_p(z)$, the superlevel set $\{\sigma_p>t\}$
is the open ball of radius $|z|$ whose measure is $\pi |z|^2$;
on the other hand,
computing $G_\newk'(s)$
from \eqref{defGk}, we can rewrite \eqref{defFp} as
\[
t=\sigma_p(z)=\frac B {c_p} G_\newk'(\pi |z|^2)^{\frac 1{p-1}}
\]
which, combined with \eqref{defmup}, yields 
$G_\newk'(\mu_p(t))=G_\newk'(\pi |z|^2)$, whence $\mu_p(t)=\pi |z|^2$
i.e. $\mu_p(t)=|\{\sigma_p>t\}|$  for every $t$ in the range of
$\sigma_p$ (that is, for every $t\in (0,B/c_p]$ since $\max \sigma_p=B/c_p$).
Finally, 
if $t>B/c_p$
we have $\{\sigma_p>t\}=\emptyset$
and in this case $\mu_p(t)=0$ by definition:
thus, summing up, we have $\mu_p(t)=|\{\sigma_p>t\}|$ for all $t>0$,
as claimed.

Now note that, when $\sigma=\sigma_p$, the two inequalities in \eqref{varb} both
become
equalities: the former, because in this case $\sigma=\sigma^*$ (hence there is
equality in \eqref{rearr} by Theorem \ref{teoF1}), and the latter
because in this case $\mu=\mu_p$ and $\mu_p$ achieves the supremum in
\eqref{varb}. Thus, the supremum in \eqref{varb} equals the right-hand
side of \eqref{rearr2}, which then follows from \eqref{varb}. Finally,
equality in \eqref{rearr2} forces a double equality in \eqref{varb}:
then necessarily $\mu=\mu_p$ (since $\mu_p$ is the only maximizer) and  
$\sigma=\sigma^*$ a.e. (by Theorem \ref{teoF1}, \eqref{rearr} being an equality).
Thus, both $\sigma^*$ and $\sigma_p$ have $\mu_p$ as distribution  function:  since
both functions are radial and radially decreasing, $\sigma^*=\sigma_p$ a.e.
and hence also $\sigma=\sigma_p$ a.e., as claimed.
\end{proof}

    \appendix
    \section{An extension of Karamata's inequality}

\begin{proof}[Proof of Proposition \ref{prop:Karamata inequality}]
We first note that $\Phi$, being convex, has a constant sign on $(0,\varepsilon)$ for some
$\varepsilon>0$, hence each of the two series in \eqref{eq:Karamata inequality}
has a well-defined sum, not necessarily finite. More precisely,
if $\Phi(0^+)$ (limit from the right) is nonzero, then the two series diverge,
both to $+\infty$ or both to $-\infty$, and \eqref{eq:Karamata inequality} is trivial
        in this case. Hence we may assume $\Phi(0^+)=0$ so that, extending $\Phi$ to $[0,1)$ by
        letting $\Phi(0):=0$, the right derivative $\Phi'(0^+)$ is well defined, with values in 
        $[-\infty,+\infty)$. Moreover, defining 
$S_0:=0$ and,
for $\newk,k\geq 1$,
\[
S_\newk:=\sum_{k=1}^\newk \left(y_k-x_k\right),\quad
c_k: = \begin{cases}
 \displaystyle   \frac{\Phi(y_k) - \Phi(x_k)}{y_k - x_k} & \text{if $x_k\not=y_k$}
 \\[3mm]
    \Phi'(x_k)\quad\text{(right derivative)} & \text{if $x_k=y_k$,}
    \end{cases}
\]
summation by parts yields
\begin{equation}
\label{sumbyparts}
    \sum_{k=1}^\newk \left( \Phi(y_k) - \Phi(x_k)\right) 
    = \sum_{k=1}^\newk c_k \left(S_k-S_{k-1}\right)=
    c_\newk S_\newk + \sum_{k=1}^{\newk-1}(c_k - c_{k+1})S_k,\quad\forall\newk\geq 1.
\end{equation}
Now two cases are possible:

\smallskip

\emph{(i) } $\Phi'(0^+)$ is finite. In this case,
$\Phi$ is Lipschitzian near $0^+$ and hence the two series in \eqref{eq:Karamata inequality}
are both convergent, due to \eqref{ass1} and the fact that $\Phi(0)=0$.
        Moreover, by the convexity of $\Phi$ we have
        $c_{k}\geq c_{k+1}\geq \Phi'(0^+)$ $\forall k\geq 1$ (in particular the sequence $c_k$ is bounded), whereas $S_{\newk}\geq 0$
        and $S_\newk\to 0$ by \eqref{ass1}. Hence, letting $\newk \to \infty$
        in \eqref{sumbyparts} leads to
        \begin{equation}\label{eq karamata_intermedia}
            \sum_{k=1}^{\infty} \left( \Phi(y_k) - \Phi(x_k)\right) = \sum_{k=1}^{\infty}(c_k-c_{k+1})S_k.
        \end{equation}
       Since as already observed $(c_k-c_{k+1})S_k \geq 0$, \eqref{eq:Karamata inequality} is proved,
       and equality occurs in \eqref{eq:Karamata inequality} if and only if
       $(c_k-c_{k+1})S_k=0$
       for every $k\geq 1$. In addition, since the interval $(0,y_1]$ is the union of 
       all the closed intervals $[x_k\wedge y_k,x_k\vee y_k]$ ($k\geq 1$), the condition
       that all the $c_k$'s are equal is equivalent to the condition that $\Phi$, restricted
       to $(0,y_1]$, is affine. If this is not the case (i.e. $c_h-c_{h+1}>0$ for at least one
        $h\geq 1$), then equality in \eqref{eq:Karamata inequality} implies 
       that $S_h=0$, which completes the proof of Proposition \ref{prop:Karamata inequality}.

\smallskip

\emph{(ii) } $\Phi'(0^+)=-\infty$. In this case, $\Phi(x)<0$ if $x>0$ is small enough,
hence the series in \eqref{eq:Karamata inequality} might each be convergent or divergent to $-\infty$. Note that each of the functions
\[
\Phi_n(t):=\max\{\Phi(t),-nt\},\quad t\in [0,1), \quad n\geq 1
\]
is convex on $[0,1)$ and such that $\Phi'_n(0^+)=-n$, hence it falls within case \emph{(i)} for
which Proposition \ref{prop:Karamata inequality} has already been proved. Thus,
\eqref{eq:Karamata inequality} with $\Phi_n$ in place of
$\Phi$ gives
\begin{equation}
\label{prev}
            \sum_{k=1}^{\infty} \Phi_n(x_k) \leq \sum_{k=1}^{\infty} \Phi_n(y_k)\quad\forall n\geq 1.
\end{equation}
Then, since $\Phi_n\searrow \Phi$ and (as observed in case \emph{(i)}) the series in \eqref{prev} are convergent,
letting $n\to\infty$ in \eqref{prev} one obtains
\eqref{eq:Karamata inequality} from monotone convergence
(the limit series in \eqref{eq:Karamata inequality}, however, need not be both
convergent).

To complete the proof, assume now
that $\Phi$ is non-affine on $(0,y_1]$ and
that equality occurs in \eqref{eq:Karamata inequality}, the two series being convergent.
Let $a\in (0,y_1)$ be any point where $\Phi$ is differentiable, chosen small enough so that
$\Phi$, restricted to the interval $[a,y_1]$, is not affine, and consider the decomposition
$\Phi(t) = \Phi_1(t) + \Phi_2(t)$, where
\begin{equation*}
            \Phi_2(t) = \begin{cases}
                \Phi(t) & \text{if $t \in [0,a]$}\\
                \Phi(a) + \Phi'(a)(t-a) & \text{if $t \in [a,1)$.}
            \end{cases}
\end{equation*}
Clearly, both $\Phi_1$ and $\Phi_2$ are convex functions on $[0,1)$, hence
the first part of Proposition \ref{prop:Karamata inequality} (which has
already been proved for any convex $\Phi$) yields
\begin{equation}
\label{prev2}
    \sum_{k=1}^{\infty} \Phi_1(x_k) \leq \sum_{k=1}^{\infty} \Phi_1(y_k),\quad
    \sum_{k=1}^{\infty} \Phi_2(x_k) \leq \sum_{k=1}^{\infty} \Phi_2(y_k).
\end{equation}
Adding these inequality yields \eqref{eq:Karamata inequality} again, but
since we are assuming that \eqref{eq:Karamata inequality} holds as an equality, 
there must be equality also in \eqref{prev2}, in particular for $\Phi_1$.
Now observe that, since $\Phi_2$ is affine on the interval $[a,1]$ but
$\Phi$ is not, on that interval $\Phi_1=\Phi-\Phi_2$ is not affine either (in particular,
$\Phi_1$ is not affine on $(0,y_1]$). But since $\Phi_1'(0^+)=0$, $\Phi_1$ falls
within case \emph{(i)}, in which Proposition \eqref{prop:Karamata inequality} has
already been proved: then, since there is equality in \eqref{prev2}, the last claim
of Proposition \eqref{prop:Karamata inequality} applies to $\Phi_1$, hence
$S_k=0$ for at least one value of $k$.
\end{proof}

\begin{remark}
If $S_k=0$ for some $k\geq 1$ and $x_{k+1}\vee y_{k+1}\leq
x_k\wedge y_k$, then equality occurs in \eqref{eq:Karamata inequality}
provided $\Phi$ is affine on the interval $(0,x_{k+1}\vee y_{k+1})$. More
elaborate examples of equality cases can be constructed, with $S_k=0$
for several non-consecutive values of $k$ and $\Phi$ piecewise affine.
In fact, with some effort, building on the previous proof it is possible to characterize
all cases where equality occurs, but we will not pursue this here.
\end{remark}
\section*{Acknowledgments}
F.N. is a fellow of the Accademia delle Scienze di Torino and a member of the Societ\`a Italiana di Scienze e Tecnologie Quantistiche (SISTEQ). 
\section*{Statements and Declarations --- Competing interests}
The authors have no financial or non-financial interests to disclose.
	

\end{document}